\input amstex
\documentstyle{amsppt}
\define\bZ{{\Bbb Z}}

\define\Cl{\text{Cliff}(TM)}
\topmatter
\title The $K$-homology class of the\\
Euler characteristic operator\\
is trivial\endtitle
\rightheadtext{The Euler characteristic operator}
\author Jonathan Rosenberg{$^\dag$}\endauthor
\leftheadtext{Jonathan Rosenberg}
\address Department of Mathematics, University of
Maryland, College Park, MD 20742, USA\endaddress
\email {\tt jmr\@math.umd.edu}\endemail
\thanks {$^\dag$}Partially supported by NSF Grant \# DMS-96-25336
and by the General Research Board of the University
of Maryland.\endthanks
\keywords $K$-homology, de~Rham operator, signature operator, Kasparov
theory
\endkeywords
\subjclass Primary 19K33; Secondary 19K35, 19K56, 58G12
\endsubjclass
\abstract On any manifold $M^n$, the de~Rham operator
$D=d+d^*$ (with respect to a complete Riemannian metric),
with the grading of forms by parity of degree, gives rise
by Kasparov theory
to a class $[D]\in KO_0(M)$, which when $M$ is closed maps
to the Euler characteristic $\chi(M)$ in $KO_0(\hbox{pt})=
\bZ$. The purpose of this note is to give a quick proof
of the (perhaps unfortunate) fact that $[D]$ is as trivial
as it could be subject to this constraint. More precisely,
if $M$ is connected, $[D]$ lies in the image of
$\bZ=KO_0(\hbox{pt})\to KO_0(M)$ (induced by the inclusion
of a basepoint into $M$).
\endabstract
\endtopmatter
\document
Let $M^n$ be a complete Riemannian manifold without
boundary (possibly compact, possibly non-compact). 
Recall that the de~Rham operator $D=d+d^*$, acting on differential
forms on $M$ (of all possible degrees) is a formally
self-adjoint elliptic operator, and that on the
Hilbert space of $L^2$ forms, it is essentially
self-adjoint \cite{Ga}. With a certain grading on the form bundle
(coming from the Hodge $*$-operator), $D$ becomes the
{\it signature operator\/}; with the more obvious grading
of forms by parity of the degree, $D$ becomes the
{\it Euler characteristic operator}. When $M$ is compact,
the kernel of $D$, the space of harmonic forms, is naturally identified 
with the real or complex\footnote{depending on what scalars one is using}
cohomology of $M$ by
the Hodge Theorem, and in this way one observes that
the index of $D$ (with respect to the parity grading) is
simply the Euler characteristic of $M$, whereas the
index with respect to the other grading is the signature
\cite{AS3}.

Now by Kasparov theory (good general references are
\cite{Bl} and \cite{Hig1}), an elliptic operator such as
$D$ gives rise to a $K$-homology class. In the case of
a compact manifold, the index of the operator is recovered
by looking at the image of this class under the map collapsing $M$ 
to a point. However, the $K$-homology
class usually carries far more information than the index
alone; for example, it determines the index of the
operator with coefficients in any vector bundle, and even determines the
families index in $K^*(X)$ of a family of twists of the operator, 
as determined by a vector bundle on $M\times X$. ($X$ here is a
parameter space.) When $M$ is non-compact, 
things are similar, except that usually there is no index, and
the class lives in an appropriate Kasparov
group $K^{-*}(C_0(M))$, which is {\it locally finite\/} $K$-homology,
{\it i.e.}, the relative group $K_*(\overline M, \{\infty\})$, 
where $\overline M$ is the one-point compactification of $M$.\footnote{Here
$C_0(M)$ denotes continuous real- or complex-valued functions on $M$
vanishing at infinity, depending on whether one is using real or complex
scalars. This algebra is {\it contra\/}variant in $M$, so a contravariant
functor of $C_0(M)$ is {\it co\/}variant in $M$.
Excision in Kasparov theory identifies $K^{-*}(C_0(M))$
with $K^{-*}(C(\overline M),\, C(\hbox{pt}))$, which is identified with
relative $K$-homology. When $\overline M$ does not have finite homotopy
type, $K$-homology here means Steenrod $K$-homology, as explained in 
\cite{KKS}.}

In the case of the signature operator, the corresponding
$K$-homology class (which lives in complex $K$-homology
in degree $n=\dim M$), is a very rich object. Rationally, it
is the Poincar\'e dual of the total $\Cal L$-class 
(the Atiyah-Singer $L$-class, which differs from the Hirzebruch
$L$-class only by certain well-understood powers of $2$),
but in addition, it also carries quite interesting integral
information \cite{KM, R, RW}. It is therefore natural
to ask if the $K$-homology class corresponding to the
Euler characteristic operator carries any similar
extra information (beyond the value of the Euler characteristic itself). 
We shall prove here that this is not the case. This result is supported
by ``popular wisdom,'' and of course one can prove a rational version of
it by using the Atiyah-Singer Theorem and characteristic class calculations,
but as far as I know the integral result has not been published before.
(However, Bruce Williams has kindly shown me a purely topological proof,
not any more elementary than the analytic one we shall give here,
based on ideas in \cite{BS, \S6}.)

We may as well work with real scalars and real $K$-theory, since
triviality of the class in real $K$-homology will imply its triviality
in complex $K$-homology as well. The fact that $D$ is odd with respect 
to the grading then gives us a class in $KO_0(M)$ (regardless of
the value of $n$). Our main result is thus:
\proclaim{Theorem} Let $M^n$ be a complete Riemannian
manifold, connected for simplicity,
let $D$ be its Euler characteristic operator,
and let $[D]\in KO_0(M)$ be its $K$-homology class
\rom{(}in locally finite $K$-homology\rom{)}. Then this class is
as trivial as one can possibly expect it to be. In other
words, if $M$ is compact, it is just
$\chi(M)\in\bZ$, with $\bZ$ embedded in $KO_0(M)$ as the
image of $KO_0(\hbox{pt})$ under inclusion of a basepoint,
and if $M$ is non-compact, the class is zero.
\endproclaim

Before getting to the proof, we need two technical results.
The first says that the $K$-homology class of the Euler characteristic
(or signature) operator depends only on the manifold and not on the choice of
a Riemannian metric.
\proclaim{Proposition} If $M^n$ is a connected manifold without boundary,
either compact or non-compact, then the class of the Euler characteristic 
operator
$[D]\in KO_0(M)$ is independent of the choice of a complete Riemannian metric
on $M$. The same statement holds for the class of the signature operator
in $K_n(M)$.
\endproclaim
\demo{Proof \rom{(}Sketch\rom{)}} 
This was proven in \cite{Hig2}, which was unfortunately
never published, though some of the ideas appear in \cite{Hig3} and \cite{Hig4}.
A proof in the more complicated context of the signature operator on 
Lipschitz manifolds appears in \cite{Hil1} (for the compact case) and
in \cite{Hil2, \S2} (for the non-compact case). For the reader's convenience
(and also because we will need a similar argument later), we summarize
the idea of Hilsum's approach.  This makes use of
the Baaj-Julg ``unbounded'' version (\cite{BJ} or \cite{Bl,
\S 17.11})  of the definition of the Kasparov $K$-homology group $KO_0(M)$. 
A class $[D]$ in $KO_0(M)$ is defined by an operator $D$ on a $\bZ/2$-graded
Hilbert space $\Cal H=\Cal H^+\oplus \Cal H^-$, where $\Cal H$ is
equipped with an action of $C_0(M)$. We shall denote 
the action by multiplication,
since in most cases of practical interest, $\Cal H$ is a space of $L^2$
sections of a vector bundle over $M$, and $C_0(M)$ acts by multiplication
operators.  The operator $D$ and the action of $C_0(M)$ are required to satisfy
the following properties:
\roster
\itemitem"{(UB1)}" The operator $D$ is self-adjoint
on $\Cal H$  (though allowed to be unbounded).
\itemitem"{(UB2)}" The operator $D$
is odd with respect to the grading of $\Cal H$, 
{\it i.e.}, $\hbox{dom}\, D =(\hbox{dom}\, D)^+\oplus (\hbox{dom}\, D)^-$,
where $(\hbox{dom}\, D)^\pm=(\hbox{dom}\, D)\cap \Cal H^\pm$, and
$D$ maps $(\hbox{dom}\, D)^+$ to $\Cal H^-$, $(\hbox{dom}\, D)^-$ to $\Cal H^+$.
\itemitem"{(UB3)}" For a dense set of $f\in C_0(M)$, $f$ maps the domain of 
$D$ into itself, and has bounded commutator with $D$.
\itemitem"{(UB4)}" For $f\in C_0(M)$, $f(1+D^2)^{-1}$ (which makes sense as a
bounded operator, by the functional calculus for unbounded self-adjoint
operators) is a compact operator.
\endroster
The class $[D]$ only depends on the homotopy class of $D$ (within the
space of all operators satisfying (UB1)--(UB4)), as made more
precise in \cite{BJ, Remarque 2.5(iv)}. In particular, if $a$ is bounded
and self-adjoint, and if $t>0$, then $D$, $tD$, and $D+a$ 
all define the same class $[D]$ \cite{BJ, Remarque 2.5(iii)}.

Given these preliminaries, the proof of the Proposition is based on the
following observations:
\roster\item The space of complete Riemannian metrics on $M$ is arcwise
connected, in fact is a convex cone (cf.\ \cite{Hil2, Corollaire 1.4}).
\item If one joins two complete Riemannian metrics $g_0$ and $g_1$ by
an arc $g_t$ of such metrics, then the corresponding Hodge $*$-operators
vary continuously, so one gets a continuous field of Hilbert spaces $\Cal H_t$
of $L^2$ forms. One can write down continuously varying explicit unitary
operators $U_t:\Cal H_0\to \Cal H_t$ that can be used to carry operators
on $\Cal H_t$ back to $\Cal H_0$.
\item If $D_t$ is the de~Rham operator on $\Cal H_t$, then it satisfies
axioms (UB1)--(UB4) (with respect to either the parity grading or,
when $n$ is even, the signature grading).
\item The operators $U_t^* D_t U_t$ satisfy the axioms for a homotopy
of unbounded Kasparov modules. Hence $[D_t]$ is independent of $t$.
\endroster
The details may be found in the references cited above.
\qed\enddemo

\proclaim{Lemma 1} Let $M$ be a connected manifold without boundary,
either compact with $\chi(M)=0$, or else non-compact.  Then
$M$ admits an everywhere non-vanishing vector field $X$ and
a complete Riemannian metric $g$ with the properties that $\Vert X\Vert=1$
and $\Vert \nabla_V X\Vert \le 1$ \rom{(}pointwise everywhere\rom{)}, for
all unit tangent vectors $V$. 
\rom{(}All norms are taken with respect to $g$, and
as usual, $\nabla$ denotes the covariant derivative with respect
to the Riemannian connection.\rom{)}
\endproclaim
\demo{Proof} The condition that either $\chi(M)=0$ or else $M$ is non-compact
guarantees that $M$ admits an everywhere non-vanishing vector field $X$.
In the compact case, choose any Riemannian metric $g'$ on $M$, and rescale
$X$ to have length 1 everywhere.  By compactness, $\Vert \nabla X\Vert' \le C$
for some constant $C>0$.  (Here $\Vert\,\cdot\,\Vert'$ denotes length with
respect to $g'$.) If we rescale $g'$ to a metric
$g$, multiplying vector lengths by $C$,
then $\frac{1}{C}X$ has unit length in the new metric.  Then if $V$ is any
unit tangent vector with respect to $g'$, $\frac{1}{C}V$ is a unit tangent
vector in the new metric $g$, and 
$$\left\Vert\nabla_{\frac{1}{C}V}(\tfrac{1}{C}X)\right\Vert
=\frac{1}{C^2}\Vert{\nabla'}_V X\Vert
=\frac{C}{C^2}\Vert{\nabla'}_V X\Vert'\le 1.$$

For the non-compact case, we use the same idea, except we need to allow more
general conformal rescaling. As in the compact case, we start with some
complete Riemannian metric $g'$ on $M$, and rescale $X$ to have length 1 
everywhere (with respect to $g'$). We then take $\langle\,\cdot\,,\,\,\cdot\,
\rangle_g=e^{2u}\langle\,\cdot\,,\,\,\cdot\, \rangle_{g'}$, where the 
function $u$ is still to be determined. Then $e^{-u}X$ will be a unit vector
field in the new metric $g$. If $V$  is a unit tangent vector with respect
to $g'$, then $e^{-u}V$ is a unit tangent vector with respect
to $g$, and we have
$$\align
\Vert \nabla_{e^{-u}V} e^{-u}X \Vert
&=e^{-u}\Vert \nabla_{V} e^{-u}X \Vert\\
&=e^{-u}\Vert \left(V\cdot e^{-u}\right)X + e^{-u}\nabla_{V}X \Vert\\
&=\Vert e^{-u}\left(-V\cdot u\right)X + e^{-u}\nabla_{V}X \Vert',\\
\intertext{(since $e^{-u}\Vert\,\,\Vert=\Vert\,\,\Vert'$)}
&=e^{-u}\Vert -\left(V\cdot u\right)X + {\nabla'}_V X
+\left(V\cdot u\right)X + \left(X\cdot u\right)V\\
&\qquad\qquad - \langle V,\, X\rangle' \hbox{grad}'(u)\Vert'\\
\intertext{(by the calculation in
\cite{LM, Proof of Theorem II.5.24, p.\ 133})}
&=e^{-u}\Vert {\nabla'}_V X
+ \left(X\cdot u\right)V - \langle V,\, X\rangle' \hbox{grad}'(u)\Vert'\\
&\le e^{-u}\left(\left\Vert{\nabla'}_V X\right\Vert'+2\left\Vert
\hbox{grad}'(u)\right\Vert'\right),\tag1
\endalign$$
with $\hbox{grad}'(u)$ here computed in the original metric $g'$. We just need
to choose $u$ so that the right-hand side of (1) is $\le 1$.

To see that this is feasible, we fix a basepoint $x_0$ and choose $u$ to
depend roughly only on the distance $r$ to $x_0$,
with respect to the metric $g'$.  (Since the distance function
may not be smooth, some smoothing may be required, but this does not affect
the basic estimates.) We choose $u$ so that 
$$e^{u(r_0)}\ge 2\max_{r\le r_0}
\left\Vert{\nabla'}_V X\right\Vert'.$$
This bounds the first term in (1) by $\frac12$. 
The second term is less of a problem,
since $e^{-u}\Vert\hbox{grad}'(u)\Vert'$ behaves like 
$$\left|{d\over dr}\left(e^{-u(r)}\right)\right|.$$
Since we may choose $u$ to increase to $+\infty$, $e^{-u}\searrow 0$,
and thus if $u$ is chosen smooth enough,
the derivative of $e^{-u(r)}$ tends to 0 as $r\to\infty$.
The new metric is complete since distances are bigger than in the old metric.
\qed\enddemo

\demo{Proof of Theorem} First observe that it is enough to prove the
theorem in the {\it non\/}-compact case. For if $M$ is
compact and connected, choose a basepoint $x_0\in M$
and let $N=M\smallsetminus\{x_0\}$. Then $KO_0(N)=
KO_0(M,\,\{x_0\})$ (recall we are using locally finite
homology), and one can show that $[D_M]$ maps
to $[D_N]$ in the relative group.  (This is not totally
trivial, as the metric on $M$ has to be rescaled on $N$
in order to give a complete metric there, but see \cite{Hig2}
and \cite{Hig3}. The point is that
the $K$-homology class in locally finite homology really only sees the
restrictions of the metric to compact sets, where all metrics are equivalent.) 
So if $[D_N]=0$ in $KO_0(N)$, $[D_M]$
must come from $KO_0(\hbox{pt})$, as claimed.

Thus we may assume $M$ is non-compact, and we will show
$[D]=0$. (The same argument would show directly that
$[D]=0$ if $M$ is compact with vanishing Euler characteristic.) 
The Proposition says we are free to choose a complete metric on $M$
as we please. So apply Lemma 1 to choose a complete metric on $M$ and a vector 
field $X$ of length 1 everywhere with respect to the given metric, so that
in addition $\Vert\nabla_V X\Vert\le 1$ for all unit tangent vectors $V$.
We identify the form bundle of
$M$ with the Clifford algebra bundle $\Cl$ of the tangent bundle, 
with its standard grading in which vector fields are sections of $\Cl^-$, and 
$D$ with the Dirac operator on $\Cl$.\footnote{Since sign conventions 
differ, we emphasize that for us, unit tangent vectors on $M$ have square $-1$
in the Clifford algebra.} (This is legitimate by \cite{LM, II, Theorem 5.12}.) 
Let $A$ be the operator
on $\Cl$ defined by {\it right\/} Clifford multiplication by $X$ on
$\Cl^+$ (the even part of $\Cl$) and by right
Clifford multiplication by $-X$ on $\Cl^-$ (the odd part). 
We use {\it right\/} Clifford multiplication since it commutes
with the symbol of $D$. Observe that with respect to the $L^2$ inner
product on sections of $\Cl$, $A$ is self-adjoint with square $1$. 
Furthermore, $A$ is odd with respect to the grading and 
commutes with multiplication by scalar-valued functions.

For $\lambda\ge 0$, let $D_\lambda=D+\lambda A$.
Since $D$ and $A$ satisfy (UB1) and (UB2), so does each $D_\lambda$. 
Since $A$ is bounded, $D_\lambda$ also satisfies (UB3) and (UB4),
and as noted above, all the $D_\lambda$ define the same Kasparov class.
\enddemo
\proclaim{Lemma 2} Let the metric on $M^n$ and the vector field $X$ satisfy
the conclusions of Lemma 1, and let $A$ and $D_\lambda$ be as above. Then
in the sense of ordering of self-adjoint operators,
$$ -n\lambda\le D_\lambda^2-(D^2+\lambda^2)\le n\lambda.$$
In particular, for $\lambda>n$, $D_\lambda^2\ge \lambda(\lambda-n)$, so
$D_\lambda$ has a bounded inverse.
\endproclaim
\demo{Proof of Lemma 2} Let $\omega$ be a section of $\Cl$, say of
$\Cl^+$. Then if $\cdot$ denotes Clifford multiplication, we have:
$$\align
D_\lambda^2\omega&=(D+\lambda A)^2\omega\\
&=(D+\lambda A)(D\omega + \lambda\omega\cdot X)\\
&=D^2\omega+\lambda D(\omega\cdot X)-\lambda (D\omega)\cdot X
-\lambda^2 \omega\cdot X\cdot X\\
&=(D^2+\lambda^2)\omega + \lambda\bigl(D(\omega\cdot X)
	-(D\omega)\cdot X \bigr).\tag2
\endalign  $$
However, from the definition of $D$ in terms of a local orthonormal
frame $e_1,\,\ldots,\,e_n$, we have:
$$\align
D(\omega\cdot X)&=\sum_j e_j\cdot\nabla_{e_j}(\omega\cdot X)\\
\intertext{which, since each $\nabla_{e_j}$ is a derivation, becomes}
&=\sum_j e_j\cdot \bigl(\nabla_{e_j}\omega\cdot X + 
	\omega\cdot\nabla_{e_j} X \bigr)\\
&=D(\omega)\cdot X + \sum_j e_j\cdot\omega\cdot \nabla_{e_j} X.
\endalign  $$
Substituting this into  equation (1), we obtain (locally)
$$D_\lambda^2\omega = (D^2+\lambda^2)\omega +\lambda\sum_j e_j
\cdot\omega\cdot\nabla_{e_j} X.$$
Taking the inner product
with $\omega$, we obtain:
$$\Vert D_\lambda\omega\Vert^2 = \Vert D\omega\Vert^2 +
\lambda^2\Vert\omega\Vert^2 + \lambda\sum_j\langle e_j\cdot
\omega\cdot\nabla_{e_j} X,\,\omega\rangle,$$
so estimating the last term, using the fact that $\Vert \nabla_{e_j} X\Vert
\le 1$, gives
$$\left\vert\left\langle \left(D_\lambda^2-D^2-\lambda^2\right)\omega,\,
	\omega\right\rangle\right\vert \le
	n\lambda \Vert\omega\Vert^2.\tag3$$
The conclusion of the Lemma follows.
A similar argument applies if $\omega$ is a section of $\Cl^-$; 
the sign of the error term in equation (2) is reversed, but equation (3)
remains unchanged.
\qed\enddemo
\demo{Proof of Theorem \rom{(}continued\rom{)}} 
We have seen that the class $[D]$ may
be defined by $D_\lambda$, for any $\lambda > 0$. In the ``bounded picture''
of Kasparov theory, the corresponding operator is 
$$\align
B_\lambda&=D_\lambda
\left(1+D_\lambda^2\right)^{-\frac12}\\
&=\frac{1}{\lambda}D_\lambda\left(\frac{1}{\lambda^2}+
\frac{1}{\lambda^2}D_\lambda^2\right)^{-\frac12}
\endalign$$
The axioms satisfied by this operator which correspond to (UB1)--(UB4)
are the following:
\roster
\itemitem"{(B1)}" It is self-adjoint, of norm $\le 1$.
\itemitem"{(B2)}" It is odd with respect to the grading of $\Cl$.
\itemitem"{(B3)}" For $f\in C_0(M)$, $fB_\lambda\sim B_\lambda f$ and 
$fB_\lambda^2\sim f$, where $\sim$ denotes equality modulo compact operators.
\endroster

We claim $B_\lambda\to A$ in the strong operator topology as 
$\lambda\to \infty$. For this we apply Lemma 2. Indeed, for any $\omega$:
$$\align
\Vert B_\lambda\omega - A\omega\Vert &=
\left\Vert \left(\frac{1}{\lambda^2}+
\frac{1}{\lambda^2}D_\lambda^2\right)^{-\frac12}
\left(A+\frac{1}{\lambda}D\right)\omega -A\omega\right\Vert \\
&\le\left\Vert \left(\frac{1}{\lambda^2}+
\frac{1}{\lambda^2}D_\lambda^2\right)^{-\frac12}\left(
\frac{1}{\lambda}D\omega\right)\right\Vert \\
&\quad +\left\Vert 
\left[\left(\frac{1}{\lambda^2}+
\frac{1}{\lambda^2}D_\lambda^2\right)^{-\frac12} - 1\right] A\omega\right\Vert.
\tag4
\endalign$$
To estimate the two terms on the right, note that by (3) we have:
$$1+\lambda^2 -C\lambda \le 1+ D_\lambda^2 \le
1+ D^2 + \lambda^2 +C\lambda .$$
In particular, as $\lambda \to \infty$,
$$\left\Vert \left(\frac{1}{\lambda^2}+
\frac{1}{\lambda^2}D_\lambda^2\right)^{-\frac12}\right\Vert
\le 1 + O\left(\frac{1}{\lambda}\right).$$ 
Since also $\frac{1}{\lambda}D\omega\to 0$, the first term on the
right in (4) goes to 0. To estimate the second term in (4),
note that on the spectral  subspace where $0\le D^2\le \mu$,
we have 
$$1+\lambda^2 -C\lambda \le 1+ D_\lambda^2 \le
1+ \mu + \lambda^2 +C\lambda,$$
and thus
$$1 - O\left(\frac{1}{\lambda}\right) \le \left(\frac{1}{\lambda^2}+
\frac{1}{\lambda^2}D_\lambda^2\right)^{-\frac12}\le
1 + O\left(\frac{1}{\lambda}\right).
$$
Hence the second term in (4) goes to $0$ for $A\omega$ in the
spectral  subspace where $0\le D^2\le \mu$. Since we can let $\mu\to\infty$,
we see $B_\lambda\to A$ in the strong operator
topology. Since $B_\lambda$ and $A$ each
satisfy the conditions (B1)--(B2), we have a candidate
for a homotopy of Kasparov modules between the Kasparov modules defined
by $D$ and by $A$. But $A=A^*=A^{-1}$ and $A$ commutes
with the action of $C_0(M)$. In other words, in the case, of $A$, we have
a degenerate module, in the sense that we can replace $\sim$ by $=$ in
condition (B3). Hence the class $[A]$ is trivial in $KO_0(M)$.

So we will be done if we can check the remaining conditions for a 
homotopy of Kasparov modules. There are only two of these. First
of all, we need to check that $1-B_\lambda^2 \to 1 - A^2=0$ in norm,
not just strongly. But $1-B_\lambda^2 = (1+ D_\lambda^2)^{-1}$, which
is bounded in norm by $(\lambda^2 - n \lambda)^{-1}$ for large
$\lambda$, by Lemma 2. Finally, we need to check that for $f\in C_0(M)$,
the commutator of $f$ and $B_\lambda$ tends to $0$ in norm. For this,
take $f$ smooth so that $[D,\,f]$ is bounded, and estimate as follows:
$$[f,\,B_\lambda] = \left[f,\,D_\lambda(1+ D_\lambda^2)^{-1/2}\right]
=[f,\,D_\lambda](1+ D_\lambda^2)^{-1/2} + D_\lambda
\left[f,\,(1+ D_\lambda^2)^{-1/2}\right].$$
Now $[f,\,D_\lambda]=[f,\,D]$ and $(1+ D_\lambda^2)^{-1/2}\to 0$ in norm,
so the first term on the right goes to $0$ in norm. As for the second
term, we have (following \cite{Bl, p.\ 199})
$$D_\lambda \left[f,\,(1+ D_\lambda^2)^{-1/2}\right]
= {1\over \pi}\int_0^\infty \mu^{-1/2}D_\lambda
\left[f,\,(1+ D_\lambda^2+\mu)^{-1}\right] \,d\mu,$$
and
$$ D_\lambda \left[f,\,(1+ D_\lambda^2+\mu)^{-1}\right] 
= D_\lambda (1+ D_\lambda^2+\mu)^{-1} \left[1+ D_\lambda^2+\mu,\,f\right]
(1+ D_\lambda^2+\mu)^{-1}.$$
Now use the fact that 
$$\left[1+ D_\lambda^2+\mu,\,f\right]=\left[D_\lambda^2,\,f\right]
=D_\lambda [D_\lambda,\,f] + [D_\lambda,\,f] D_\lambda = D_\lambda [D,\,f]
+[D,\,f] D_\lambda.$$
We obtain that
$$\multline
D_\lambda \left[f,\,(1+ D_\lambda^2+\mu)^{-1}\right] \\
={D_\lambda^2\over 1+D_\lambda^2+\mu} [D,\,f] {1\over 1+D_\lambda^2+\mu}
+ {D_\lambda\over 1+D_\lambda^2+\mu} [D,\,f] {D_\lambda\over 1+D_\lambda^2+\mu},
\endmultline$$
which can be bounded in norm by $2\left\Vert[D,\,f]\right\Vert$  times
$(1+\lambda^2-n\lambda +\mu)^{-1}$, for large $\lambda$. But
$${1\over \pi}\int_0^\infty \mu^{-1/2}{1\over 1+\lambda^2-n\lambda +\mu} \,d\mu
= {1\over \sqrt{1+\lambda^2-n\lambda }} \to 0\quad\hbox{as }\lambda\to 0.$$
So $\left\Vert[f,\,B_\lambda]\right\Vert\to 0$ and we have a homotopy of Kasparov
modules. Thus $[D]=[A]=0$.
\qed\enddemo

\refstyle{A}

\enddocument